\documentclass[10pt]{article}
\usepackage{amsmath,amsthm,amsfonts}

\input epsf.tex
\newdimen\epsfxsize
\newdimen\epsfysize
\def\qed{\vrule height5pt width3pt depth.5pt}

\theoremstyle{plain}
\newtheorem{thm}{Theorem}[section]
\newtheorem{cor}[thm]{Corollary}

\newtheorem{prop}[thm]{Proposition}

\newtheorem{rem}{Remark}[section]

\begin{document}

\title{Invariants from the Linking Number}

\author{H. A.  Dye   }

\maketitle

\begin{abstract}
We explore a family of invariants obtained from linking numbers. This is a family of Kauffman finite type invariants. 
\end{abstract}

\section{Introduction}
\begin{figure}[htb] \epsfysize = 0.5 in
\centerline{\epsffile{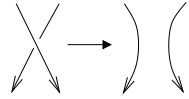}}
\caption{Vertical Smoothing}
\label{fig:smooth}
\end{figure}

Given a virtual knot $K$ smoothing a crossing vertically (as shown in Figure \ref{fig:smooth}) results in a 2 component link $L$. 
We compute the linking number of the link $L$ by summing the sign of the crossings that involve edges from two different components. For classical knots, 
the linking number is always divisible by two. In the virtual case, this does not occur and we obtain linking numbers that are congruent to both zero and one mod two. This allows us to compute a new type of invariant that can be computed directly from either the Gauss code or signed chord diagrams, and is a Kauffman finite type invariant of virtual links and knots. This invariant is related to the invariants of Henrich \cite{allison} and Manturov \cite{m1}, \cite{m2}. 

\begin{figure}[htb] \epsfysize = 1 in
\centerline{\epsffile{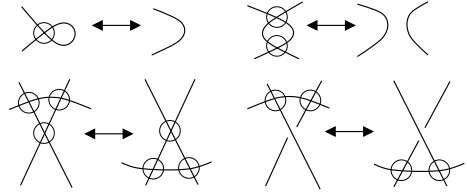}}
\caption{Classical Reidemeister move}
\label{fig:moves}
\end{figure}

\begin{figure}[htb] \epsfysize = 1 in
\centerline{\epsffile{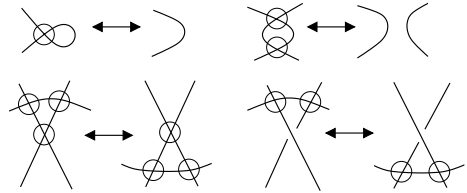}}
\caption{Virtual Reidemeister move}
\label{fig:vmoves}
\end{figure}

Recall that a virtual link diagram is a decorated immersion of $n$ copies of $S^1$ into the plane that contains two types of crossings, classical crossings (indicated by over-under markings) and virtual crossings (indicated by a solid encircled X). Two virtual link diagrams are said to be equivalent if they are related by a sequence of classical and virtual Reidemeister moves. These moves are illustrated in Figures \ref{fig:moves} and \ref{fig:vmoves}. In this paper, virtual knots and virtual knot diagrams will simply be referred to as knot and knot diagrams for convenience.

\begin{figure}[htb] \epsfysize = 0.5 in
\centerline{\epsffile{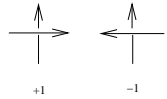}}
\caption{Crossing sign}
\label{fig:sgn}
\end{figure}

Let $K$ denote a virtual knot diagram and let $c$ denote a crossing of $K$.  We denote the sign of the crossing as $sgn(c)$ and the sign of the crossing is evaluated as shown in Figure \ref{fig:sgn}.
The writhe of a virtual knot $K$, $w(K)$, is defined to be:
\begin{equation*}
w(k) = \sum_{c \in K} sgn(c).
\end{equation*}

Fix an orientation of the virtual knot $K$ and a choose a crossing $c$.
Let
$K_c$ denote the two component link diagram  (with an inherited orientation) obtained by smoothing the crossing $c$ vertically.  Let $K_i$ denote the $i^{th}$ component of $K_c$ and let $K_1 \cap K_2 $ denote the set of crossings in $K_c$ that contain one strand from each component. We define
 $ L(K_c) $ to be a linking number of the two components
where
\begin{equation*}
L(K_c) = \sum_{c \in K_1 \cap K_2 }  sgn(K_c).
\end{equation*}
Further, let $ \bar{L} (K_c) = L(K_c)$ mod 2. 

\begin{rem} The classical definition of linking number either divides $ \bar{L_c} $ by two or sums over crossings where $K_1$ is the overcrossing and $K_2$ is the undercrossing or vice versa. \end{rem}

We define $ \gamma (K) $ to be a sum over all crossings in $K$:
\begin{equation} \label{definvar}
\gamma (K) = \sum_{c} t^{\bar{L} (K_C)} sgn(c)
\end{equation}
Note that for any virtual knot diagram $K$, $ \gamma (K) $ has the form $ a + b t $ where $a, b \in 
\mathbb{Z}$. That is, $ \gamma (K) $ is an element of a free left $ \mathbb{Z}$ module over the set $ {1, t }$.

We define $ \bar{ \gamma} (K) $:
\begin{equation}
\bar{ \gamma} (K)  = \gamma (K) mod 2.
\end{equation}

\begin{figure}[htb] \epsfysize = 0.5 in
\centerline{\epsffile{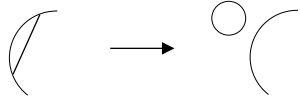}}
\caption{Reidemeister I move}
\label{fig:r1chord}
\end{figure}

\begin{thm}[Invariance] The sums $ \gamma(K) $ and $\bar { \gamma} (K) $ are invariant under all classical and virtual Reidemeister moves except Reidemeister move I.\end{thm}
\textbf{Proof:} We use chord diagrams to show the invariance of $ \gamma (K)$. An introduction to chord diagrams can be found in \cite{gpv}.
We note that $ \gamma (K) $ is not invariant under the Reidemeister I move. The introduction of a single Reidemeister I type twist is indicated on a chord diagram by an isolated chord. Smoothing along an isolated chord produces an unlinked two component link diagram, contributing $ \pm 1 $ to the sum, as shown in Figure \ref{fig:r1chord}.

A Reidemeister II move in a knot diagram introduces two crossings with opposite sign.  
\begin{figure}[htb] \epsfysize = 1.5 in
\centerline{\epsffile{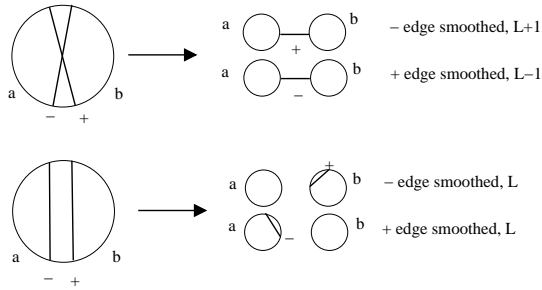}}
\caption{Reidemeister II move}
\label{fig:r2chord}
\end{figure}

The two possible chord diagrams obtained after the introduction of a Reidemeister II move are shown in Figure \ref{fig:r2chord}. To obtain the chord diagram corresponding to the removal of a Reidemeister II move, simply remove the pair of chords.
In the first type of Reidemeister II move (shown at the top of Figure 
\ref{fig:r2chord}) smoothing vertically along the negative edge results in a link ($K_{n}$) with linking number $L+1$ and smoothing along the positive edge results in a link ($K_{+}$) with linking number $L-1 $.  Hence
\begin{equation*}
\bar{L} (K_{n}) = \bar{L} (K_{+}).
\end{equation*}
Now
\begin{equation*}
\gamma (K) = \sum_{c} t^{\bar{L} (K_c) } sgn(c) \end{equation*}
and
\begin{equation*}
\gamma (K) = t^{\bar{L} (K_+)} - t^{\bar{L} (K_{n})} + \sum_{c \neq +,n } t^{\bar{L} (K_c) } sgn(c).
\end{equation*}

 In the second type of Reidemeister II move (shown at the bottom of Figure 
\ref{fig:r2chord}), smoothing along the positive and negative edges produces links with linking number $L$.  Because the edges have opposite sign, these contributions cancel out in the sum $ \gamma (K)$.

\begin{figure}[htb] \epsfysize = 2 in
\centerline{\epsffile{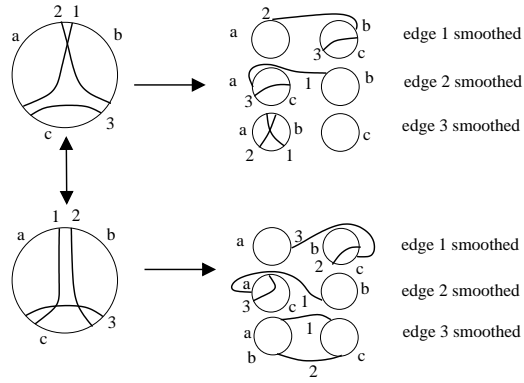}}
\caption{Reidemeister III move, type I}
\label{fig:r3chordt1}
\end{figure}

\begin{figure}[htb] \epsfysize = 2 in
\centerline{\epsffile{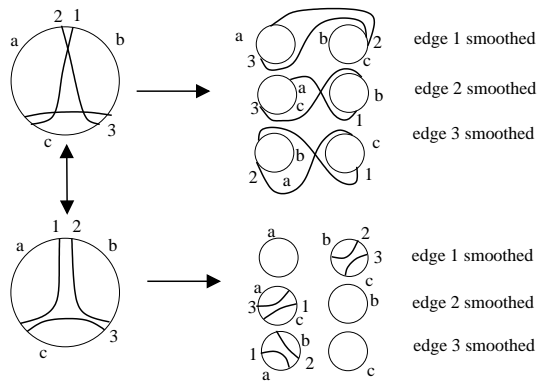}}
\caption{Reidemeister III move, type 2}
\label{fig:r3chordt2}
\end{figure}

There are two sets of chord diagrams of the Reidemeister III moves, based on the directionality of the edges. We will refer to these sets as type 1 and type 2, Figures \ref{fig:r3chordt1} and \ref{fig:r3chordt2} respectively. In these diagrams, numbers are used to distinguish individual chords while letters indicate  (possible) arrangements of endpoints of chords. On the left hand side of Figure \ref{fig:r3chordt1}, we show the chord diagrams corresponding to both sides of a Reidemeister III move, type 1. On the right hand side, we show schematics of the chords diagrams of the links obtained by vertically smoothing each of the three crossings involved in a Reidemeister III move. 

For a vertical smoothing of edge $i$, that the linking number of the resultant links (although the links themselves may not be equivalent) are equivalent modulo two. As a result, both result in the same contribution to $ \gamma (K)$. 
The same argument applies to the second type of Reidemeister III move as shown in 
Figure \ref{fig:r3chordt2}. 
As a result, $ \gamma (K) $ and $\bar{ \gamma } (K) $ are both invariant under the classical and virtual Reidemeister moves.\qed

We immediately realize the following propositions.

\begin{prop} For a virtual knot $K$, $ \gamma (K) = a + (2k) t$ where $a$ and $k$ are integers.
\end{prop}
\textbf{Proof:} Let $K$ be a virtual knot diagram and $C_K$ denote a chord diagram corresponding to $K$. The invariant $ \gamma (K)$ counts, up to sign, the number of chords that are intersected by an even number of chords and the number of chords that are intersected by an odd number of chords. The simplest chord diagram contains zero chords and has $0$ even chord and $0$ odd chords. A diagram with $1$ chord contains $1$ even chord and $0$ odd chords. We assume that a diagram with $n$ chords has an even number of odd chords, denoted by $k_o$. Let $k_e $ denote the number of even chords. Suppose that we introduce a single new chord to the diagram. 
We have four cases based on the parity of the number of even and odd chords intersected.

 \textbf{Case 1:} The new chord intersects both an even number of odd and even chords, denoted $2o$ and $2e$.  We note that the new chord is an even chord. Now:
\begin{gather*}
k_o \rightarrow k_o + 2e - 2o \\
k_e \rightarrow k_e - 2e + 2o+1
\end{gather*}
As result, the total number of odd chords has changed by an even number.

\textbf{Case 2:} The new chord intersects an even number of odd chords and an odd number even chords, denoted $2o$ and $2e+1$.  We note that the new chord is an odd chord. Now:
\begin{gather*}
k_o \rightarrow k_o - 2o + (2e+1) + 1 \\
k_e \rightarrow k_e + 2o -(2e + 1)
\end{gather*}
The total number of odd chords has changed by an even number. 

 \textbf{Case 3:} The new chord intersects an odd number of odd chords and an even number of even chords. This also changes the total number of odd chords by an even number. The computation is analogous.

  \textbf{Case 4:} The new chord intersects an odd number of odd chords and an odd number of even chords. The total number of chords changes by an even number.\qed

\begin{prop} For a virtual knot $K$, $ \gamma (K)_{t=1} $ is the writhe. In particular, for a classical knot diagram, $ \gamma (K) = w(K)$. \end{prop}
\textbf{Proof:} Clear. \qed

\begin{prop} \label{fintype} Let $K$ denote a virtual knot diagram with a positive crossing $l$. Let $K'$ denote a virtual knot diagram where the crossing $l$ is switched to a negative crossing.  Then $ \gamma (K) - \gamma (K') = \pm 2 t^{\bar{L} (K_l) } $ and $ \bar{ \gamma }(K)  - \bar{ \gamma} (K') =0 $. \end{prop}
\textbf{Proof:}
Let $K$ denote a knot diagram with a positive crossing $l$. Let $K'$ denote a knot diagram where the crossing $l$ is switched to a negative crossing. For any crossing $c$ except $l$, the smoothed diagrams $K_c$ and $K'_c $ will differ only at the crossing $l$. If the crossing $l$ does not involve both componnents then it is clear that $ \bar{L} (K_c) =
\bar{L} (K'_c) $. However, if the crossing $l$ does involve both components of the smoothing link diagram then $ L (K_c) = m+1 $ and $L (K'_c) = m-1 $ for some $m \in 
\mathbb{Z}$. Hence, $ \bar{L} (K_c) = \bar{L} (K'_c) $. If the crossing $l$ is smoothed then $ \bar{L} (K_l) = \bar{L} (K'_L) $. 

Then

\begin{gather} \label{f1}
\gamma (K) = \sum_{c \in K} t^{\bar{L}(K_c)} sgn(c) \\
= t^{ \bar{L} (K_L) }+ \sum_{c \in K, c \neq l} t^{ \bar{L} (K_c)} sgn (c).
\end{gather}

Similarly:
\begin{gather} \label{f2}
\gamma (K') = \sum_{c \in K'} t^{\bar{L}(K'_c)} sgn(c)  \\
= -t^{\bar{L} (K'_L)} + \sum_{c \in K', c \neq l} t^{ \bar{L} (K'_c)} sgn (c).
\end{gather}
Combining equations \ref{f1} and \ref{f2}, we obtain
\begin{equation*}
\gamma (K) - \gamma (K') = 2 t^{\bar{L} (K_L) } 
\end{equation*}
and
\begin{equation*}
\bar{ \gamma } (K) - \bar{ \gamma (K') } =0. \qed
\end{equation*}

\begin{cor} Let $K$ denote a virtual knot diagram with a positive crossing $l$. Let $K'$ denote a virtual knot diagram where the crossing $l$ is switched to a negative crossing.  Then $ \gamma (K) + \gamma (K') \equiv 0 $ modulo 2. 
\end{cor}
Let $S$ denote a subset of the classical crossings in the virtual knot diagram $D$. Let $D_S$ denote the diagram obtained by switching the set of crossings in $S$. Let $|S|$ denote the cardinality of $S$.
Recall that a finite type invariant $v$ has degree $n$ if for $ |S| =n $ then 
$ v(D) = (-1)^|S| v(D_s) $ \cite{kvirt}.
By this definition, proposition \ref{fintype} shows that $ \bar{ \gamma} (K) $ is a Kauffman finite type invariant of degree one \cite{gpv} \cite{allison}.

\section{Extending the invariant}
We can extend $\bar{\gamma} $ to form a sequence of invariants by smoothing pairs of crossings and summing the knot diagrams obtained in this manner.
A chord has odd parity (denoted by a $\bullet $ in diagrams) if an odd number of chords intersect the chord. Otherwise, a chord has even parity (denoted by a $ \circ $ in diagrams). Parity is invariant under the classical and virtual Reidemeister moves. (Note that this is the parity defined by Manturov in \cite{m2}.)  Let $p$ denote a pair of intersecting chords with opposite parity and let $P$ denote the set of all interesecting chords with opposite parity. The virtual knot diagram obtained by smoothing the crossings in the diagram that correspond to the pair of chords is $K_p$.
We define $ \varsigma(K) $, a formal sum of diagrams obtained by smoothing the original diagram along selected pairs of crossings. We select pairs of crossings that correspond to intersecting chords with opposite parity.
\begin{equation}
\varsigma(K) = \sum_{p \in P}  K_p.
\end{equation}
We now define $ \bar{\gamma}_2 (K) =  \gamma (t^2 \varsigma(K) ) \text{mod 2}$. 
\begin{rem}For each diagram $K_p$, we obtain a polynomial of the form $ a + b t$ for $a,b \in \mathbb{Z}$. The evaluation of $ \gamma (t^2 \varsigma(K) ) $ is $ \sum_{p \in P} t^2 \gamma (K_p)$. This becomes $ \sum_{p \in P} t^2 (a_p + b_p t) $. We consider this sum mod 2 in order to obtain invariance under the Reidemeister moves. Unfortunately, since $ b_p $ is always even, $ \bar{ \gamma}_2 (K) $ is either $0$ or $ t^3$.  \end{rem}

We prove that $ \bar{ \gamma}_2 $  is invariant under the classical and virtual Reidemeister moves.

\begin{thm} For a virtual knot diagram $K$, $\bar{ \gamma }_2 (K) $ is invariant under the Reidemeister moves. 
\end{thm}
\textbf{Proof:}
Let $K$ be a virtual knot diagram and let $C(K)$ represent the corresponding chord diagram. In the chord diagrams for this proof, we observe the following conventions. Ennumerated, signed edges represent the chords corresponding to the crossings involved in the studied moves. The labeled edge $g$ represents a generic edge that intersects the ennumerated edges. The letters $a,b,c,d,e$ represent sequences of chord endpoints and are used to calculate and compare $\gamma (K)$ from different chord diagrams.

Note that  a crossing introduced by a Reidemeister I move has even parity. The edge is isolated and does not intersect any edges. Hence, it does not contribute to $\varsigma(K)$. 
(Notice that if $K$ is the unknot with a single twist, $ \varsigma (K) $ is the empty sum and $ \bar{ \gamma}_2 (K) = 0 $.)

\begin{figure}[htb] \epsfysize = 2.5 in
\centerline{\epsffile{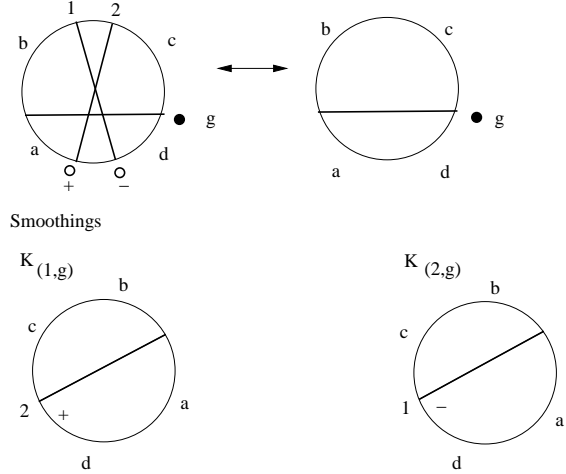}}
\caption{Reidemeister II move}
\label{fig:g2r2t1}
\end{figure}

We now consider the Reidemeister II move. A Reidemeister II move introduces two oppositely signed crossings. In a chord diagram, these crossings are represented as a pair of edges with the same parity. Without loss of generality we will assume that the crossings have even parity as shown in Figure \ref{fig:g2r2t1}.
In the diagram, the edges $1$ and $2$ represent the crossings involved in the Reidemeister II. An edge $g$ with opposite parity will either intersect both $1$ and $2$ or intersect either. The smoothing pairs $(1,g)$ and $(2,g)$ both contribute to the sum $\varsigma(K)$. Note that $ \bar{ \gamma} (K_{ (1,g)}) = \bar{ \gamma } (K_{ (2,g) }) $ for a net contribution of zero to $ \bar{ \gamma}_2 (K) $.

For the second type of Reidemeister II move, we obtain an analogous pair of diagrams.

We consider the Reidemeister III move. There are two sets of chord diagrams corresponding to the Reidemeister III move and two cases for the parity of the chords involved. 
We first consider the case where all three chords in the Reidemeister III move has the same parity. Without loss of generality, we let all three chords have even parity. In the computation of $ \varsigma (K) $, we do not obtain any diagrams from pairs of chords in the Reidemeister III move. A contribution to $ \varsigma (K) $ would involve one chord from the Reidemeister III move and a chord of different parity that intersects the Reidmeister III move chords.

\begin{figure}[htb] \epsfysize = 3 in
\centerline{\epsffile{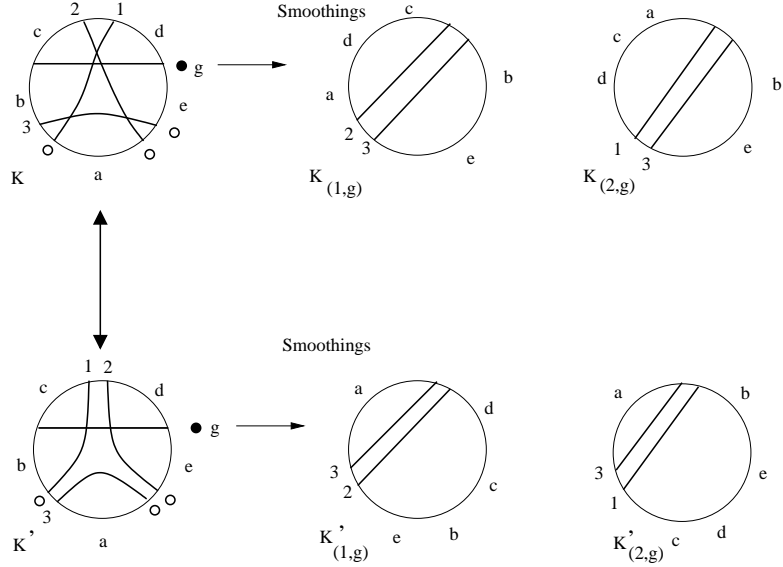}}
\caption{Reidemeister III move, type 1}
\label{fig:g2r3t1}
\end{figure}

In Figure \ref{fig:g2r3t1}, we examine a type 1 Reidmeister move where all three chords have even parity. We consider the chord 
$g$ as shown in Figure \ref{fig:g2r3t1}. The chord diarams obtained by smoothing the pairs $(1,g)$ and $(2,g) $ contribute to $ \varsigma $. We observe that
$ \gamma (K_{(i,g)} ) = \gamma ( K' _ {(i,g)} )$.
We observe that in the case of other smoothing pairs, the Reideister III move remains intact.

\begin{figure}[htb] \epsfysize = 3 in
\centerline{\epsffile{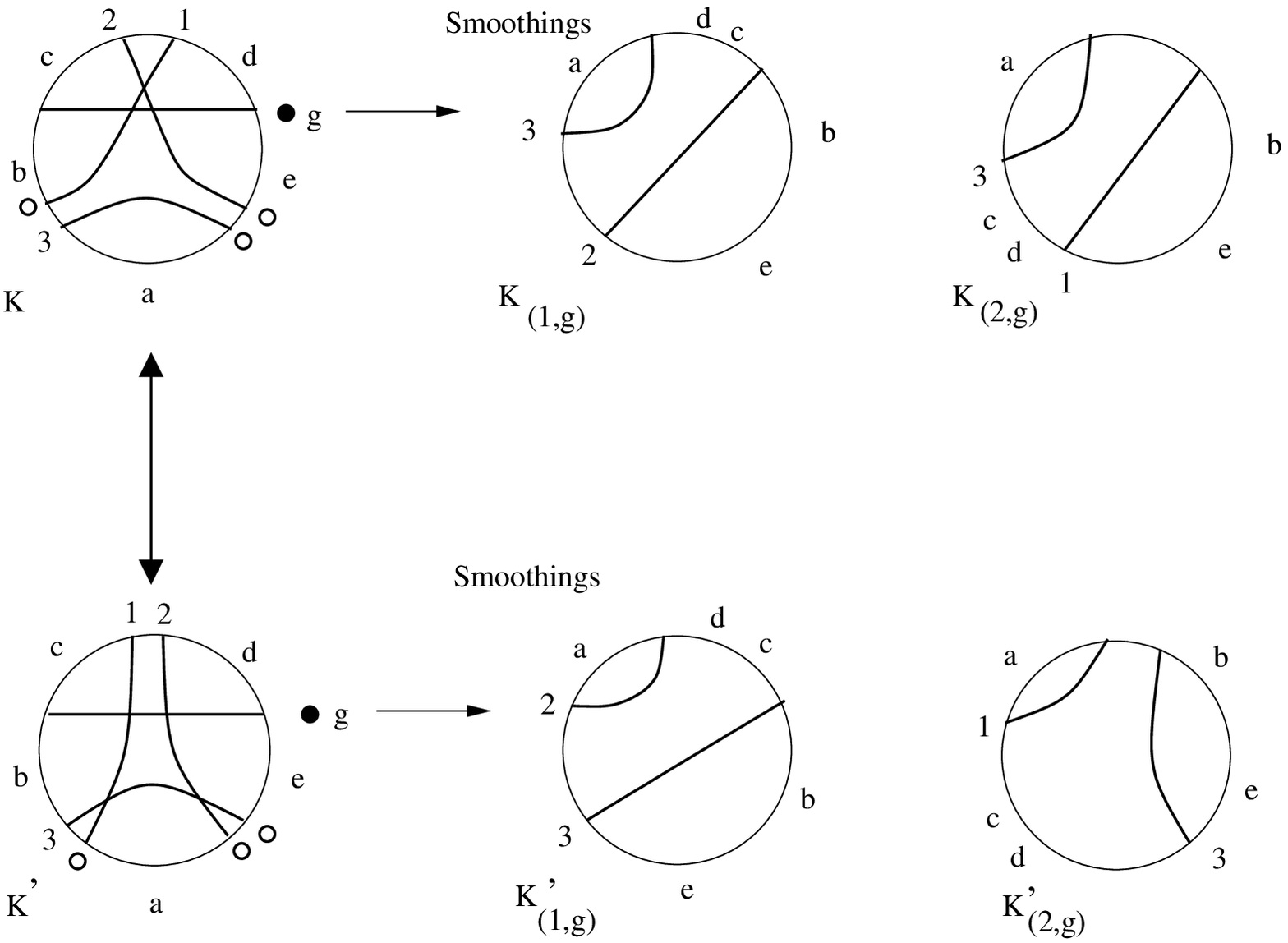}}
\caption{Reidemeister III move, type 2}
\label{fig:g2r3t2}
\end{figure}

\begin{figure}[htb] \epsfysize = 3 in
\centerline{\epsffile{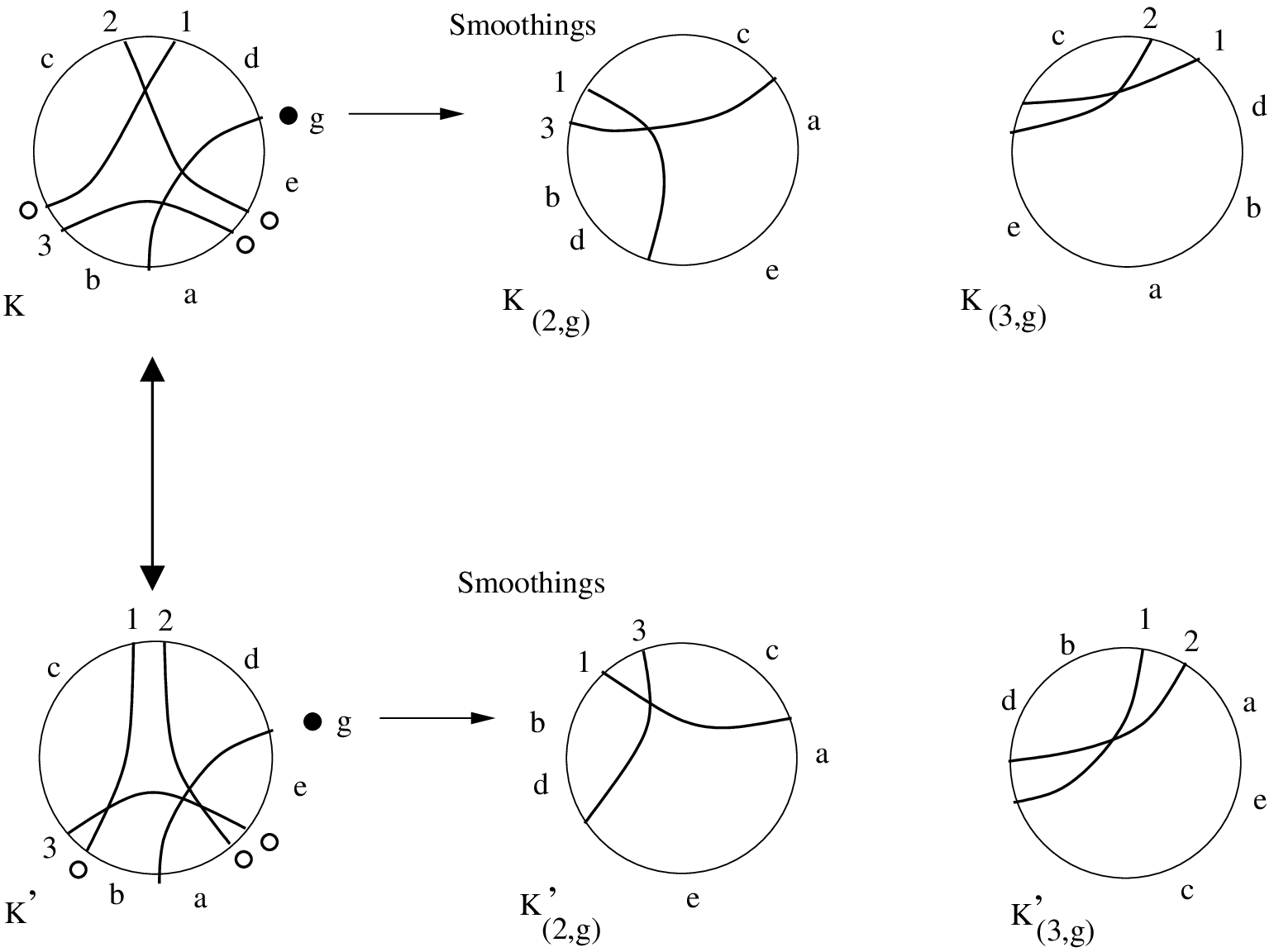}}
\caption{Reidemeister III move, type 2, case 2}
\label{fig:g2r3t3}
\end{figure}

In Figures \ref{fig:g2r3t2} and \ref{fig:g2r3t3}, we consider the Reidemeister move III, type 2 where all chords have even parity. There are two potential ways for a chord $g$ with opposite parity to intersect the Reidemeister III move.

\begin{figure}[htb] \epsfysize = 3 in
\centerline{\epsffile{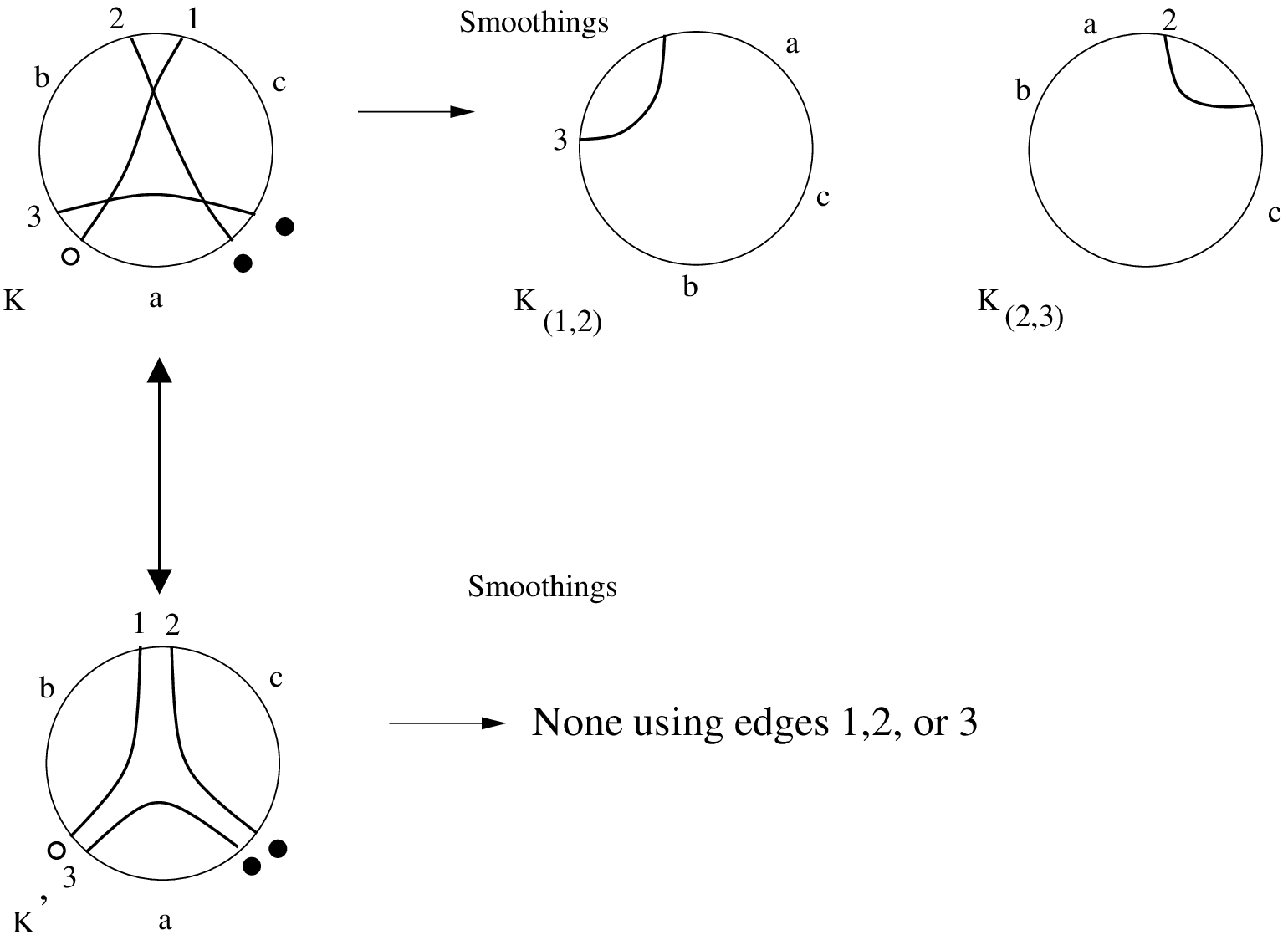}}
\caption{Reidemeister III move, type 1, edges with different parity}
\label{fig:diffparityg2r3t1}
\end{figure}

\begin{figure}[htb] \epsfysize = 3 in
\centerline{\epsffile{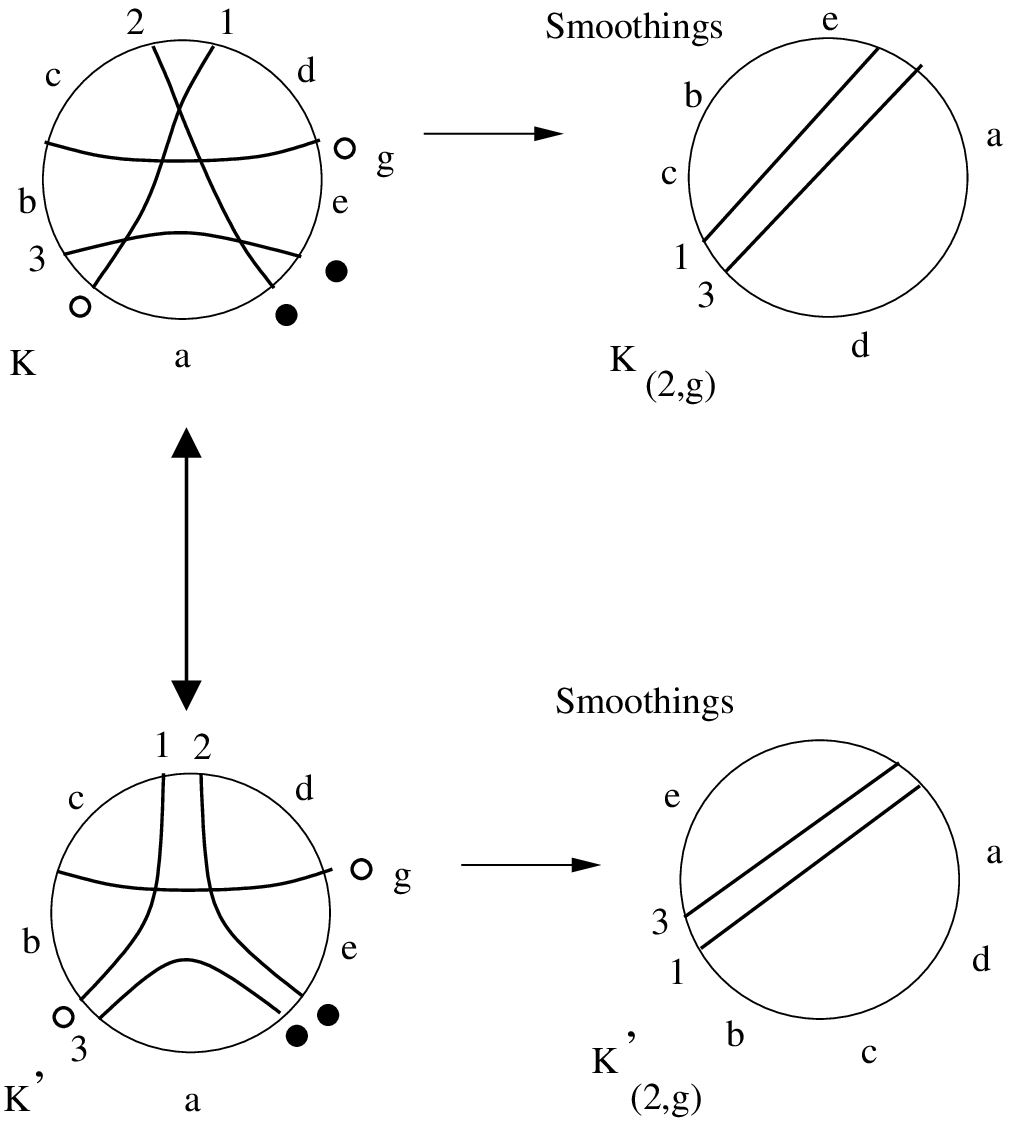}}
\caption{Reidemeister III move, type2, edges with different parity}
\label{fig:diffparityg2r3t2}
\end{figure}

\begin{figure}[htb] \epsfysize = 3 in
\centerline{\epsffile{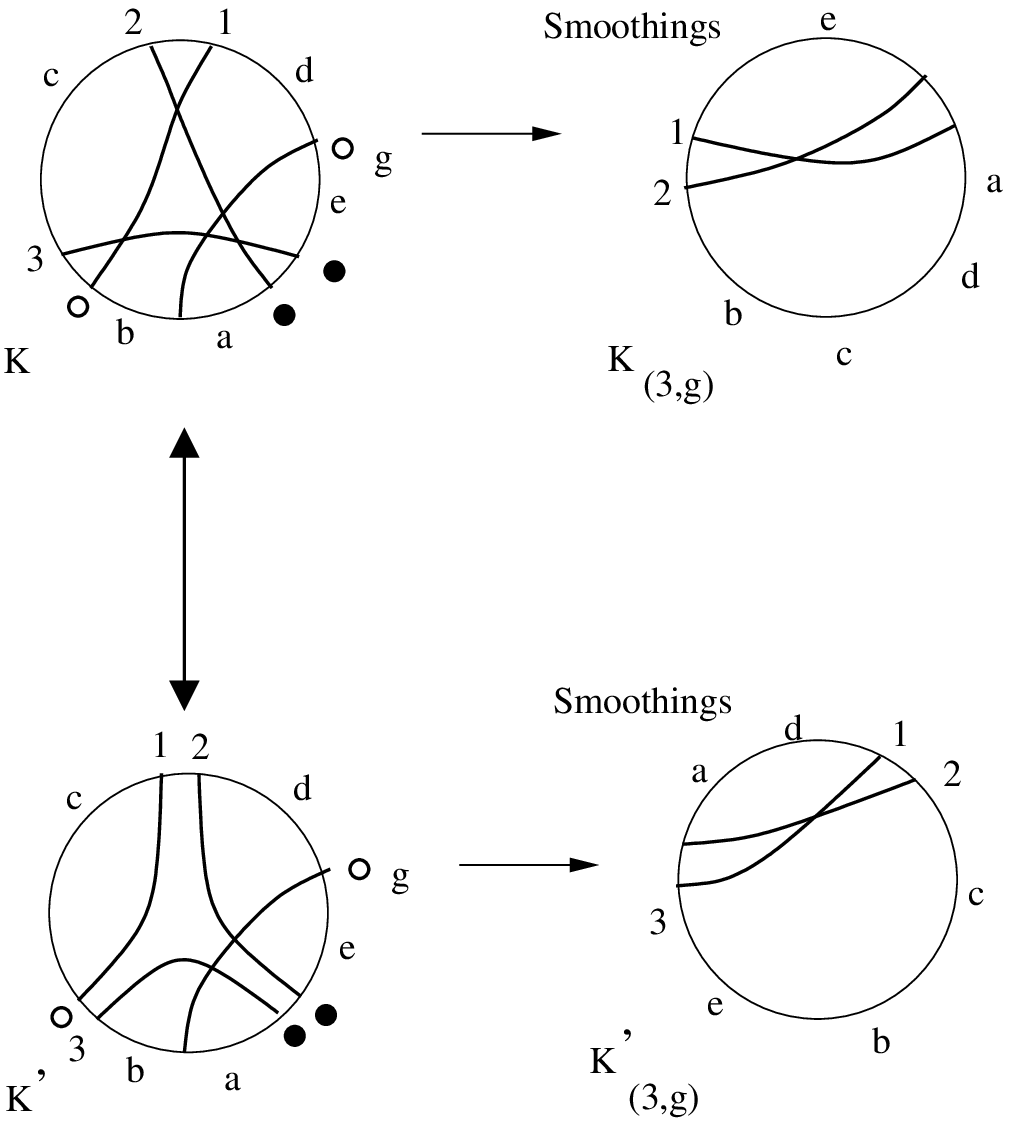}}
\caption{Reidemeister III move, type 2, case 2, edges with different parity}
\label{fig:diffparityg2r3t3}
\end{figure}

In Figures \ref{fig:diffparityg2r3t1}, \ref{fig:diffparityg2r3t2}, and \ref{fig:diffparityg2r3t3}, we consider the Reidemeister III move where the chords have differing parity. The most interesting case is shown in 
 Figure \ref{fig:diffparityg2r3t1}. From $K$, we obtain two diagrams $K_{(1,2)} $ and $K_{ (2,3)} $ which differ only by a Reidemeister I move. Note that $ \bar{ \gamma} (K_{ (1,2)} ) = 
\bar{ \gamma} (K_{(2,3)}) $, and as a result these terms cancel out in $ \gamma_2 $. From $K'$, we obtain no diagrams resulting in the same net contribution to $ \gamma_2 $.

In Figure \ref{fig:diffparityg2r3t2} and \ref{fig:diffparityg2r3t3}, we examine the cases obtained from a type 2 Reidemeister III moves where the edges have different parity.

We observe that the net contributions are equivalent and $ \bar{ \gamma}_2 $ is invariant under the classical and virtual Reidemeister moves.\qed
\begin{prop} The invariant $  \bar{ \gamma}_2 $ is a finite type invariant of degree one. \end{prop}
\textbf{Proof:} We consider the knot diagrams $K$ and $K'$ which are related by switching exactly one crossing, say $s$. To compute $\varsigma$, we expand pairs of crossings with opposite parity. If these pairs contain $s$, we denote the diagrams obtained from $K$ and $K'$ as $K_{p_i} $ and $K'_{p_i} $ respectively. Note that $K_{p_i} $ and $K'_{p_i} $ are the same diagram, so that $ \gamma (K_{p_i}) = \gamma (K'_{p_i}) $. If the pair does not contain $s$, we denote the diagrams obtained from $K$ and $K'$ as $K_{q_i} $ and $K'_{q_i} $ respectively. Note that $K_{q_i} $ and $K'_{q_i} $ differ by exactly one crossing.
Now by Proposition \ref{fintype}, $ \gamma (K_{p_i}) - \gamma (K'_{p_i}) \equiv 0$ modulo two. Hence $ \bar{ \gamma}_2 (K) - \bar{ \gamma}_2 (K') = 0 $ and $ \bar{ \gamma}_2$ is finite type invariant of degree one.\qed

\section{Examples}
We compute a variety of examples of the invariants $ \gamma $ and $ \gamma_n $ in this section.
\subsection{Trefoil}
\begin{figure}[htb] \epsfysize = 0.5 in
\centerline{\epsffile{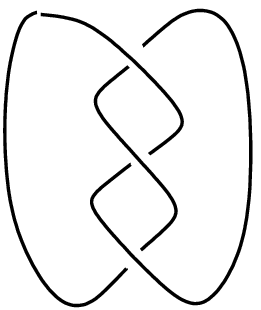}}
\caption{Trefoil}
\label{fig:tref}
\end{figure}

The trefoil shown in Figure \ref{fig:tref} has writhe $-3$. 
We observe that $ \gamma (K)= -3 $ which is equivalent to the writhe.

\subsection{Virtual Trefoil}

\begin{figure}[htb] \epsfysize = 0.5 in
\centerline{\epsffile{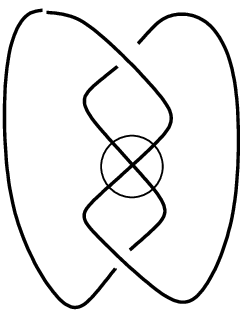}}
\caption{Virtual Trefoil}
\label{fig:vtref}
\end{figure}

The virtual trefoil, with writhe $-2$ is shown in Figure \ref{fig:vtref}.
The value of the $ \gamma $ invariant is $-2t$ and $ \gamma (K) |_{t=1}  = -2 $.

\subsection{Miyazawa's Knot}
\begin{figure}[htb] \epsfysize = 0.5 in
\centerline{\epsffile{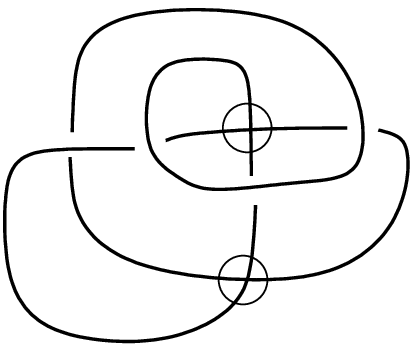}}
\caption{Miyazawa's Knot}
\label{fig:miya}
\end{figure}
In Figure \ref{fig:miya} is Miyazawa's knot and $ \gamma (K) = 2 + 2t $.
\subsection{A knot with $ \gamma (K) = 2k+3 + (2n+2)t$}
\begin{figure}[htb] \epsfysize = 0.75 in
\centerline{\epsffile{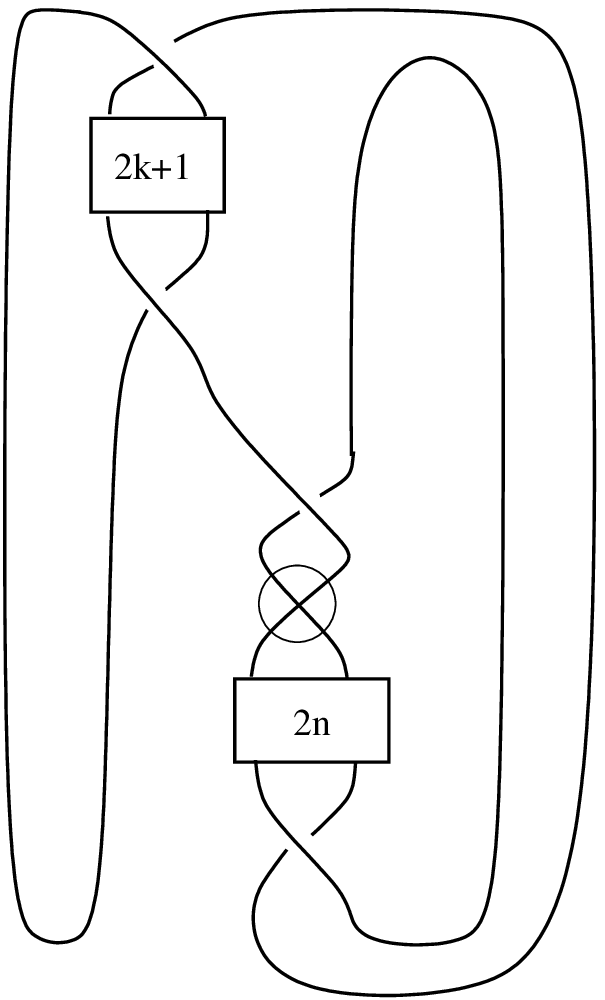}}
\caption{Knot with $\gamma (K) = (2k+3) + (2n+2) t $}
\label{fig:arbknot}
\end{figure}
For the diagram shown in Figure \ref{fig:arbknot}, complete the labeled tangles with $2k+1 $ and $2n$ positive twists, respectively. This produces a virtual knot diagram with 
$\gamma (K) = (2k+3) + (2n+2) t $ for $n,k>0 $. 

\subsection{Computing $\gamma_2 $}
\begin{figure}[htb] \epsfysize = 0.75 in
\centerline{\epsffile{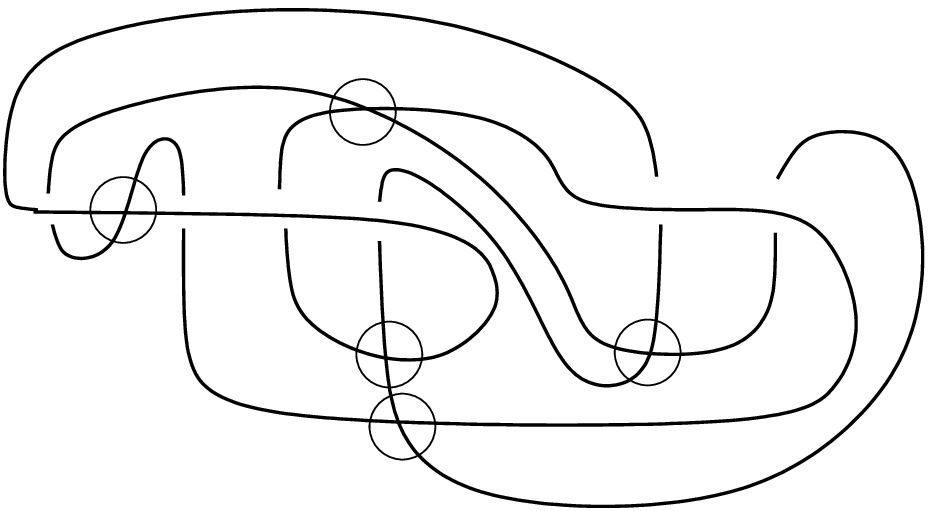}}
\caption{A knot with non-zero $ \gamma_2 $.}
\label{fig:gamma2knot}
\end{figure}
For the knot shown in Figure \ref{fig:gamma2knot}, the $ \gamma $ invariant is $ 2 + 4t $ and the $ \bar{ \gamma}_2 $ is $t^2$.

\section{Conclusion}

The invariants $ \gamma $ and $ \bar{\gamma}_2 $ can be extended by taking additional characteristics of the chord diagram into consideration. We
can extend the invariants by either enhancing the parity or considering formal sums of diagrams obtained by smoothing pairs of crossings. 
To ehance parity, we could consider over/under markings that are usually indicated on the chord diagrams by the presence of arrowheads at one end of the chord. Counting the chord intersections using the arrowheads to mark orientations will augment $ L (K_c) $. 
In the second case, it may be possible to extend $\gamma $ by smoothing particular quadruples of crossings. These ideas will be investigated in future work.

\end{document}